\documentclass[10 pt]{article}  

\usepackage{graphics} 
\usepackage{epsfig} 
\usepackage{mathptmx} 
\usepackage{times} 

\usepackage[margin=1in]{geometry}
\usepackage{cite}
\usepackage{amsmath}
\usepackage{amssymb}
\usepackage{amsfonts}
\usepackage[utf8]{inputenc}
\usepackage[english]{babel}
\usepackage{xcolor}
\usepackage{comment}
\usepackage{bbold}
\usepackage{amsthm}

\newtheorem{theorem}{Theorem}

\newtheorem{lemma}[theorem]{Lemma}

\newtheorem{assumption}{Assumption}
\newtheorem{remark}{Remark}
\usepackage{ulem}

\title{\LARGE \bf
Inner approximations of the maximal positively invariant set for polynomial dynamical systems$^1$
}
\date{April 20, 2019}

\begin{document}

\author{Antoine Oustry$^2$, Matteo Tacchi$^3$, and Didier Henrion$^4$}
\footnotetext[1]{The research was funded by the French company R\'{e}seaux de Transport d'\'{E}lectricit\'{e} (RTE).}
\footnotetext[2]{Antoine Oustry is with \'{E}cole Polytechnique, Route de Saclay, 91128 Palaiseau, France; he is also with R\'{e}seaux de Transport d'\'{E}lectricit\'{e}, Immeuble WINDOW - 7C, Place du Dôme 92073 Paris La Défense, France.
        {\tt\small antoine.oustry@polytechnique.edu}}%
\footnotetext[3]{Matteo Tacchi is with the company R\'{e}seaux de Transport d'\'{E}lectricit\'{e}, Immeuble WINDOW - 7C, Place du D\^ome 92073 Paris La D\'efense, France; he is also with LAAS-CNRS, 7 avenue du colonel Roche, 31400 Toulouse, France.
 {\tt\small tacchi@laas.fr}}%
 \footnotetext[4]{Didier Henrion is with the LAAS-CNRS, 7 avenue du colonel Roche, 31400 Toulouse, France; he is also with the Faculty of Electrical Engineering, Czech Technical University in Prague, Technick\'a 2, 166 27 Prague, Czechia.
 {\tt\small henrion@laas.fr}}%

\maketitle
\thispagestyle{empty}
\pagestyle{empty}

\begin{abstract}
The Lasserre or moment-sum-of-square hierarchy of linear matrix inequality relaxations is used to compute inner approximations of the maximal positively invariant set for continuous-time dynamical systems with polynomial vector fields. Convergence in volume of the hierarchy is proved under a technical growth condition on the average exit time of trajectories. Our contribution is to deal with inner approximations in infinite time, while former work with volume convergence guarantees proposed either outer approximations of the maximal positively invariant set or inner approximations of the region of attraction in finite time.
\end{abstract}

\section{Introduction}

This paper is an effort along a research line initiated in \cite{hk14} for developing convex optimization techniques to approximate sets relevant to non-linear control systems subject to non-linear constraints, with rigorous proofs of convergence in volume. The approximations are obtained by solving numerically a hierarchy of semidefinite programming or linear matrix inequality (LMI) relaxations, as proposed originally by Lasserre in the context of polynomial optimization \cite{l10}. Convergence proofs are achieved by exploiting duality between non-negative continuous functions and Borel measures, approximated respectively with sums of squares (SOS) of polynomials and moments, justifying the terminology moment-SOS or Lasserre hierarchy. In the context of control systems, the primal moment formulation builds upon the notion of occupation measures \cite{lhpt08} and the dual SOS formulation can be classified under Hamilton-Jacobi techniques \cite{bcht17}.

Previous works along this line include inner approximations of the region of attraction \cite{khj13}, outer approximations of the maximal positively invariant (MPI) set \cite{khj14}, as well as outer approximations of the reachability set \cite{ghmt16}. 
These techniques were applied e.g. in robotics \cite{mvtt14} and biological systems \cite{rfmfd18}.
In \cite{hk14,khj13} the regions of attraction are defined for a finite time horizon, which is a technical convenient framework since the occupation measures have then finite mass. To cope with an infinite time horizon and MPI sets, a discount factor was added in \cite{khj14} so that the mass of the occupation measure decreases fast enough when time increases. In \cite{ghmt16}, the mass was controlled by enforcing a growth condition on the volume of complement sets. This condition, difficult to check a priori, can be validated a posteriori using duality theory.

It must be emphasized here that, in general, the infinite time hoziron setup is more convenient for the classical Lyapunov framework and asymptotic stability, see e.g. \cite{c11} and references therein, whereas the finite time horizon setup is more convenient for approaches based on occupation measures. In the current paper, we make efforts to adapt the occupation measure framework to an infinite time horizon setup, at the price of technical difficulties similar to the ones already encountered in \cite{ghmt16}. Contrary to the outer approximations derived in \cite{khj14}, we have not been able to use discounted occupation measures for constructing inner approximations. Instead, the technical device on which we relied is a growth condition on the average exit time of trajectories.

The main contributions of this work are:

\begin{enumerate}
	\item A {moment-SOS} hierarchy for constructing inner approximations of the MPI set  for a polynomial dynamical system with semialgebraic constraints;
    \item A {detailed, self-contained,} rigorous proof of convergence of the hierarchy, under an assumption on the average exit time of trajectories.
\end{enumerate}

Section~\ref{sec:Statement} presents the problem statement. Section~\ref{sec:Inner} describes the MPI set inner approximation method with proof of convergence under appropriate assumptions. Numerical results are analyzed in Section~\ref{sec:Numerics}. Conclusion and future work are discussed in Section~\ref{sec:Conclusions}.

\section{Problem statement}\label{sec:Statement}

Consider the autonomous system
\begin{equation}
   \dot{x}(t) = f(x), \quad  x \in \text{int}(X) \subset \mathbb{R}^n, \quad t \in [0,+\infty [
   \label{ODE}
\end{equation}
with a given polynomial vector field $f$ of degree $d_0$. The state trajectory $x(.)$ is constrained to the interior $\text{int}(X)$ of a nonempty compact basic semi-algebraic set
$$X:= \{x \in \mathbb{R}^n : g_i(x) \geq 0, \: i = 1,...,n_X \}$$
where the $g_i$ are given polynomials of degree $d_i$. We define $\partial X := X \setminus \text{int}(X)$.

The  vector  field $f$ is  polynomial  and  therefore  Lipschitz  on  the  compact  set $X$. As  a
result, for any $x_0 \in \text{int}(X)$, there exists a unique maximal solution $x(.|x_0)$ to ordinary differential equation \eqref{ODE} with initial condition $x(0|x_0) = x_0$. The time interval on which this solution is defined contains the time interval on which $x(.|x_0)\in \text{int}(X)$. 

For any $t \in \mathbb{R}_+ \cup \{+\infty\}$, we define the following set:
$$X_t := \{x_0 \in \text{int}(X): \forall s \in [0,t],\: x(s|x_0) \in \text{int}(X) \}$$
With this notation, $X_\infty$ is the set of all initial states generating trajectories staying in $\text{int}(X)$ \textit{ad infinitum}: $X_\infty$ is the MPI set included in $\text{int}(X)$. Indeed, for any $x_0 \in X_\infty$ and $t\geq 0$, by definition, $x(t\vert x_0) \in X_\infty$. {We make the following assumption implying that $X_\infty$ has non-empty interior:

\begin{assumption} $X$ contains a Lyapunov-stable equilibrium point for $f$. \end{assumption}
}

The complementary set $X_t^c := \text{int}(X) \setminus X_t$ is the set of initial conditions generating trajectories reaching the target set $\partial X$ at any time before $t$: this is the region of attraction of $\partial X$ with free final time lower than $t$. The complementary set $X^c_{\infty}$ is the region of attraction of  $\partial X$ with free and unbounded final time. 

In this paper we want to approximate the MPI set $X_\infty$ from inside as closely as possible.

\section{Inner approximations of maximal positively invariant set}\label{sec:Inner}

{This section presents an infinite dimensional linear programming (LP) problem and a hierarchy of convex linear matrix inequality (LMI) relaxations yielding a converging sequence (in the sense of the Lebesgue measure) of inner approximations of the MPI set.}

However, due to the infinite time horizon, such a strong result is available only under some assumptions. It is based on the primal formulation of the MPI set computation problem.

For a given $x_0 \in X^c_{\infty}$, we define the exit time as
$$ \tau(x_0): = \inf \{t \geq 0 : x(t|x_0) \notin \text{int}(X)\}.$$
In the rest of this paper we make the assumption that the average exit time of trajectories leaving $\text{int}(X)$ is finite:
\begin{assumption}
$\overline{\tau} := \frac{1}{\lambda(X)}\int_{X^c_{\infty}} \tau(x)dx < + \infty$.
\label{Assumptionfinite}
\end{assumption}

{
\begin{remark}
This assumption is necessary for the rigorous proof of convergence of the sequence of approximations of $X_\infty$. It is difficult to check a priori.
We will show however that, independently of this assumption, the validity of our approximations can be checked numerically a posteriori. 
\end{remark}
}

\subsection{Primal LP}

For a given $T \in \mathbb{R}_+$, we define the following infinite-dimensional LP 
\begin{equation}
\begin{array}{rcl}
p^T \:= & \sup\limits_{\substack{\mu_0,\hat{\mu_0},\mu \in \mathcal{M}^+(X) \\ \mu_\partial \in \mathcal{M}^+(\partial X)}} & \mu_0(X) \vspace*{-1.5em}\\
&& \text{div}(f\mu)+\mu_\partial=\mu_0  \\
&\text{s.t.}&  \mu_0+\hat{\mu}_0 = \lambda \\
&& \mu(X) \leq T\: \lambda(X)
\end{array}
\label{eqn:Pa}
\end{equation}
with $\mathcal{M}^+(A)$ denoting the cone of non-negative {elements of the vector space $\mathcal{M}(A)$ of Borel measures supported on the set $A$. With the notation $\mathcal{M}(A)=\mathcal{C}(A)'$ we emphasize that the vector space $\mathcal{M}(A)$ is the topological dual to the vector space $\mathcal{C}(A)$ of continuous functions on $A$. The duality of $v \in \mathcal{C}(A)$ with  $\mu \in \mathcal{M}(A)$ is the integration $\langle \mu,v \rangle := \int_A v(x)d\mu(x)$.}
\begin{remark}
Here, $T$ is introduced to ensure that all the feasible measures have a finite norm in total variation \mbox{$\Vert \mu \Vert_{TV} := \mu(X) < +\infty$}. Otherwise, the optimization problem would be ill-posed.
\end{remark}

{Note that problem \eqref{eqn:Pa} is linear in the decision variables which are the four measures.}
The two following lemmas link the infinite-dimensional LP \eqref{eqn:Pa} and the MPI set $X_\infty$.

\begin{lemma}
Assuming that $T \geq \overline{\tau}$, we have $p^T \geq \lambda(X^c_\infty)$.
\label{LemmaInf}
\end{lemma}
\begin{proof}
\begin{itemize}
    \item $\mu_0 := \lambda_{X^c_\infty}$, $\hat{\mu}_0 := \lambda - \mu_0 = \lambda_{X_\infty}$
    \item $\mu := A \mapsto \int_{X^c_\infty} \int_0^{\tau(x_0)} \mathbb{1}_A(x(t\vert x_0)) \; dt \; dx_0$
    \item $\mu_\partial =: A \mapsto \int_{X^c_\infty} \mathbb{1}_A(x(\tau(x_0)\vert x_0))dx_0$
\end{itemize}
define a feasible quadruplet. Indeed, one has :
\begin{itemize}
\item $\mu(X) = \int_{X^c_\infty} \left(\int_0^{\tau(x_0)}dt \right) dx_0 = \overline{\tau} \: \lambda(X) \leq T \lambda(X) $
\item the first constraint in \eqref{eqn:Pa} is satisfied, since $\forall v \in \mathcal{C}^1(X)$,
\end{itemize}
\begin{align*}
    \langle \text{div}(f\mu), v \rangle & = - \int_{X^c_\infty} \int_0^{\tau(x_0)} \nabla v(x(t\vert x_0)) \cdot f(x(t\vert x_0)) dt dx_0 \\
    & = - \int_{X^c_\infty} \left( v(x(\tau(x_0)\vert x_0)) - v(x_0) \right) dx_0 \\
    & = \langle \mu_0 - \mu_\partial, v \rangle.
\end{align*}
then, $p^T \geq \mu_0(X) = \lambda_{X^c_\infty}(X) = \lambda(X^c_\infty)$.
\end{proof}
\begin{lemma}
For any quadruplet $(\mu_{0},\hat{\mu}_{0},\mu_,\mu_\partial)$ feasible in \eqref{eqn:Pa}, $\mu_0$ is supported on $X^c_\infty$, i.e. $\mu_0(X_\infty) = 0$.
\label{LemmaSup}
\end{lemma}

The proof of this lemma uses the following assumption on the MPI set:
\begin{assumption} \label{assum:enter}
$\forall x \in {\partial X_\infty} \cap \partial X,\, f(x) \cdot n(x) < 0$, where $n(x)$ stands for the unit normal vector to $\partial X$ pointing towards $\mathbb{R}^n\setminus X$.

In words, at all points where ${\partial X_\infty}$ is tangent to ${\partial X}$, the trajectories strictly enter $X$. Up to the choice of $X$, this assumption is reasonable for any physical system.
\end{assumption}
\begin{proof}
Let $(\mu_0,\hat{\mu}_0,\mu,\mu_\partial)$ be a feasible quadruplet for \eqref{eqn:Pa}. Let $\nu := \text{div}(f\mu) = \mu_0 - \mu_\partial \in \mathcal{M}(X)$. For $x \in \mathbb{R}^n$, let
$$ \varphi(x) := \begin{cases} K \exp\left(-\frac{1}{1 - \vert x \vert^2}\right) \qquad \text{if } \vert x \vert < 1 \\ 0 \qquad \qquad \qquad \qquad  \text{ else} \end{cases} $$
where $K > 0$ is such that $\int \varphi \; d\lambda = 1$. Then, for $\epsilon > 0$ and $x \in \mathbb{R}^n$, let:
\begin{itemize}
    \item $\varphi_\epsilon(x) := \frac{1}{\epsilon} \, \varphi\left(\frac{x}{\epsilon}\right) \geq 0$
    \item $\mu_\epsilon(x) := \displaystyle\int_{X} \varphi_\epsilon(y - x) \; d\mu(y) \geq 0$
    \item $\nu_\epsilon(x) := \text{div}(f\mu_\epsilon)(x)$.
\end{itemize}
According to the theory of mollifiers, $\varphi$, $\varphi_\epsilon$, $\mu_\epsilon$ and $\nu_\epsilon$ are smooth compactly supported functions, and for any $w \in \mathcal{C}^0(\mathbb{R}^n)$ compactly supported,
$$ \int_{\mathbb{R}^n} w(x) \, \mu_\epsilon(x) \; dx \underset{\epsilon \to 0}{\longrightarrow} \int_{X} w(x) \; d\mu(x)$$
from which it directly follows that for $v \in \mathcal{C}^1(\mathbb{R}^n)$ compactly supported
\begin{align*}
    \int_{\mathbb{R}^n} v(x) \, \nu_\epsilon(x) \; dx & = \int_{\mathbb{R}^n} v(x) \, \text{div}(f\mu_\epsilon)(x) \; dx \\
    & = - \int_{\mathbb{R}^n} \nabla v(x) \cdot f(x) \, \mu_\epsilon(x) \; dx \\
    & \underset{\epsilon \to 0}{\longrightarrow} - \int_{\mathbb{R}^n} \nabla v(x) \cdot f(x) \; d\mu(x) \\
    & = \int_{\mathbb{R}^n} v(x) \; d\nu(x).
\end{align*}
{ By density of $\mathcal{C}^1_c(\mathbb{R}^n)$ in $\mathcal{C}^0_ c(\mathbb{R})$ with respect to the supremum norm
$\Vert .\Vert_{L^\infty(\mathbb{R}^n)}$
}, this implies that $\nu_\epsilon \lambda$ weakly converges (in the sense of measures) to $\nu$.

For a given $\delta > 0$ consider the set
$$
X_\delta := \left\{ x \in X_\infty : \inf\limits_{y \in \partial X} \vert x - y \vert > \delta \right\}.
$$
By definition, $X_\delta \cap \partial X = \emptyset$, and then for any Borel set $A \subset X_\delta$, one has $\nu(A) = \mu_0(A)$. In particular, $\nu(\partial X_\delta) = \mu_0(\partial X_\delta) = 0$ since $\mu_0 \leq \lambda$. Then, we can apply the Portmanteau lemma (equality marked with a $\ast$) to $\nu(X_\delta)$:
\begin{align*}
    \mu_0(X_\delta) & = \nu(X_\delta) \\
    & \stackrel{\ast}{=} \lim\limits_{\epsilon \to 0} \int_{X_\delta} \nu_\epsilon(x) \; dx \\
    & \stackrel{\text{def}}{=} \lim\limits_{\epsilon \to 0} \int_{X_\delta} \text{div}(f \mu_\epsilon)(x) \; dx \\
    & = \lim\limits_{\epsilon \to 0} \int_{\partial X_\delta} f(x) \cdot n_\delta(x) \, \mu_\epsilon(x) \; dx
\end{align*}
where $n_\delta$ stands for the unit normal vector to $\partial X_\delta$ pointing towards $X_\delta^c$, according to Stokes' theorem. Now, let $\Delta$ be the function
$$ \begin{array}{ccl}
    \partial X_\infty \cap \partial X & \longrightarrow & \mathbb{R}_+ \\     x & \longmapsto & \sup \left\{\begin{array}{c}
        \Delta > 0 , \forall \delta \in (0,\Delta), \forall y \in \partial X_\delta \\
        \vert x - y \vert < \Delta \Longrightarrow f(y) \cdot n_\delta(y) \leq 0
     \end{array}\right\}.
\end{array} $$
In words, $\Delta(x)$ is the largest range around $x$ in which the $f \cdot n_\delta$ are non-positive.
According to Assumption \ref{assum:enter}, $f$ being continuous, $\Delta$ takes only positive values. Moreover, due to the regularity of $f$, $X$ and $X_\infty$, $\Delta$ is continuous on the compact set $\partial X_\infty \cap \partial X$, therefore it attains a minimum $\Delta^\ast > 0$.

Let $\delta \in (0,\Delta^\ast)$, $x \in \partial X_\delta$. Then, there are two possibilities:
\begin{itemize}
    \item either $x \in \partial X_\infty$, and then by positive invariance of $X_\infty$, $f(x)\cdot n_\delta(x) \leq 0$;
    \item or $\inf\limits_{y \in \partial X} \vert x - y \vert = \delta < \Delta^\ast$, and by definition of $\Delta^\ast$, $f(x)\cdot n_\delta(x) \leq 0$.
\end{itemize}
It follows that for any $x \in \partial X_\delta$, $f(x) \cdot n(x) \leq 0$. Thus, one obtains
$$ \int_{\partial X_\delta} f(x) \cdot n_\delta(x) \, \mu_\epsilon(x) \; dx \leq 0 $$
and after letting $\epsilon$ tend to $0$, we have $ \mu_0(X_\delta) \leq 0 $, which means, by non-negativity of $\mu_0$, that $\mu_0(X_\delta) = 0$.

{Eventually, since $X_\delta \subset X_\infty$ and $\mu_0 \leq \lambda$, one has
\begin{align*}
\mu_0(X_\infty) & = \mu_0(X_\infty) - \mu_0(X_\delta) \\
&= \mu_0(X_\infty \setminus X_\delta) \\
& \leq \lambda(X_\infty \setminus X_\delta) \\
&= \lambda(\{x \in X_\infty \; ; \, \inf\limits_{y \in \partial X} \vert x - y \vert \leq \delta \}) \\
& \underset{\delta \to 0}{\longrightarrow} \lambda(X_\infty \cap \partial X) = 0
\end{align*}
which leads to the conclusion that $\mu_0(X_\infty) = 0.$}
\end{proof}

\begin{theorem}
Assuming that $T \geq \overline{\tau}$, the infinite-dimensional LP \eqref{eqn:Pa} has a value $p^T = \lambda(X^c_\infty)$. Moreover the supremum is attained, and the $\mu_0$ component of any solution is necessarily the measure $\lambda_{X^c_\infty}$.
\end{theorem}
\begin{proof}
This is a straightforward consequence of lemmas \ref{LemmaInf} and \ref{LemmaSup}.
\end{proof}

\subsection{Dual LP}

For a given $T \in \mathbb{R}_+$, { we cast problem \eqref{eqn:Pa} as particular instance of a primal LP in the canonical form:
\begin{equation}
\begin{array}{rcl}
p^T \:= & \sup\limits_{\phi \in K_1} & \langle \phi,c \rangle_1\\
&\text{s.t.}&  - \mathcal{A}\phi + b \in K_2 \\
\end{array}
\label{eqn:PaBis}
\end{equation}
with
\begin{itemize}
\item the vector space $E_1 := \mathcal{M}(X)^3 \times \mathcal{M}(\partial X)$ and its cone $K_1 := E_1^+$ of non-negative elements;
\item the vector space $F_1 := \mathcal{C}^0(X)^3 \times \mathcal{C}^0(\partial X)$ and $L_1 := F_1^+$;
\item the duality $\langle .,. \rangle_1 : E_1 \times F_1 \rightarrow \mathbb{R}$, given by the integration of continuous functions against Borel measures, since $E_1$ is the dual of $F_1$;
\item the decision variable $\phi:=(\mu_0, \hat{\mu}_0, \mu, \mu_\partial) \in E_1$ and the objective function $c:=(1,0,0,0) \in F_1$;
\item $E_2 := \mathbb{R} \times \mathcal{C}^1(X)' \times \mathcal{M}(X)$, $K_2 := \mathbb{R}_+ \times \{0\} \times \{0\} \subset E_2$ and the right hand side vector $b:=(T \lambda(X),0, \lambda)$;
\item $F_2 := \mathbb{R} \times \mathcal{C}^1(X) \times \mathcal{C}^0(X)$, $L_2 := \mathbb{R}_+ \times \mathcal{C}^1(X) \times \mathcal{C}^0(X)$;
\item the linear operator $\mathcal{A}: E_1 \rightarrow E_2$ given by
$$\mathcal{A}(\mu_0, \hat{\mu}_0, \mu, \mu_\partial) := \begin{pmatrix}\mu(X) \\
\text{div}(f\mu)+\mu_\partial - \mu_0\\
\mu_0+\hat{\mu}_0 \\
\end{pmatrix}. $$
\end{itemize}

Note that both spaces $F_1,F_2$ are equipped with the weak topologies $\sigma(F_1,E_1)$, $\sigma(F_2,E_2)$ and the spaces $E_1,E_2$ are equipped with the weak-* topologies $\sigma(E_1,F_1)$, $\sigma(E_2,F_2)$. Using the same notations, the dual of the primal LP \eqref{eqn:PaBis} in the canonical form reads:
\begin{equation}	
\begin{array}{rcl}
d^T \:= & \inf\limits_{\psi \in L_2} & \langle b,\psi \rangle_2\\
& \text{s.t.} & \mathcal{A}'\psi - c \in L_1  \\
\end{array}
\label{eqn:DaBis}
\end{equation}

with 

\begin{itemize}
\item the Lagrange dual variable $\psi:=(u,v,w) \in F_2$;
\item the adjoint linear operator $\mathcal{A}':F_2 \rightarrow F_1$ given by:
$$\mathcal{A}'(u,v,w) := \begin{pmatrix}
w - v \\
w \\
u - \nabla v \cdot f \\
v\vert_{\partial X} \\
\end{pmatrix}. $$
\end{itemize}

Using our original notations, the dual LP of problem \eqref{eqn:Pa} then reads:}
\begin{equation}
\begin{array}{rcl}
d^T \:= & \inf\limits_{\substack{u \in \mathbb{R}_+ \\ v \in \mathcal{C}^1(X) \\ w \in \mathcal{C}^0(X)}} & \displaystyle \int_{X} \left(w(x) + u \: T\right) \: d\lambda(x) \vspace*{-1.5em}\\
&& \nabla v \cdot f (x) \leq u, \:\forall x \in X \\
& \text{s.t.} &  w(x) \geq v(x)+1, \forall x \in X \\
&& w(x) \geq 0, \forall x \in X \\
&& v(x) \geq 0, \forall x \in \partial X.
\end{array}
\label{eqn:Da}
\end{equation}

{\begin{lemma}
Let $(0,v,w)$ be a feasible triplet for problem \eqref{eqn:Da}. Then, the set $\hat{X}_\infty := \{x \in \text{int}(X) : v(x) < 0\}$ is a positively invariant subset of $X_\infty$.
\label{lem_pos_inv}
\end{lemma}
\begin{proof}
Since $X_\infty$ is the MPI set included in $X$ and $\hat{X}_\infty \subset X$ by definition, it is sufficient to prove that $\hat{X}_\infty$ is positively invariant.

Let $x_0 \in \hat{X}_\infty$. Then, for any $t > 0$, it holds
$$v(x(t\vert x_0)) =  v(x_0) + \int_0^t \nabla v \cdot f(x(s\vert x_0)) \ ds  \leq v(x_0) < 0 $$
using constraint $\nabla v \cdot f (x) \leq u = 0$.

We still have to show that $x(t\vert x_0)$ remains in $\text{int}(X)$ at all times $t \geq 0$. If not, then there exists a ${t_\partial} > 0$ such that $x({t_\partial} \vert x_0) \in \partial X$ according to the intermediate value theorem, the trajectory being of course continuous in time. However, by feasibility of $(0,v,w)$, one then has $v(x({t_\partial}\vert x_0)) \geq 0$, which is in contradiction with the fact that $v(x(t \vert x_0)) < 0$ for all $t > 0$ which we just proved.

Thus, we obtain that for all $t > 0$, $x(t\vert x_0) \in \text{int}(X)$ and $v(x(t\vert x_0)) < 0$, i.e. $x(t\vert x_0) \in \hat{X}_\infty$. 
\end{proof}}
{
\begin{remark}
For a feasible triplet $(u,v,w)$, if $u \neq 0$, then there is no guarantee that the solution of \eqref{eqn:Da} yields an inner approximation of $X_\infty$. However, it still gives access to inner approximations of the $X_t$, $t \in \mathbb{R}_+$, and we will show that under Assumption \ref{Assumptionfinite}, these approximations converge to $X_\infty$.
\end{remark}}

\begin{lemma}
For any triplet $(u,v,w)$ feasible in \eqref{eqn:Da}, for any $ t>0$, $\hat{X}_t := \{x_0 \in \text{int}(X), v(x_0) + ut < 0 \} \subset X_t$.
\label{lem_inner_approx}
\end{lemma}
\begin{proof}
Let $(u,v,w)$ be a feasible triplet in \eqref{eqn:Da} and let $x_0$ be a element of $X_t^c$ for a given $t>0$. 

By definition of $X_t$ we know that $t\geq \tau(x_0)$, where $\tau$ is the exit time, and that for any $s \in [0,\tau(x_0)],x(s|x_0) \in X$. Thanks to the first constraint in \eqref{eqn:Da}, we can therefore say that for any $s \in [0,\tau(x_0)], (\nabla v\cdot f) (x(s|x_0)) \leq u$. Hence for any $s \in [0,\tau(x_0)], v(x(s|x_0)) \leq v(x_0) + us$. In particular, we deduce that $$v(x(\tau(x_0)|x_0))\leq v(x_0) + u\tau(x_0) \leq v(x_0) + ut$$
As $x(\tau(x_0)|x_0) \in \partial X$, we know that $v(x(\tau(x_0)|x_0)) \geq 0$ and thus $v(x_0) \geq - ut$. This proves that $$X_t^c \subset \{x_0 \in \text{int}(X), v(x_0) \geq - ut \}$$ hence $\hat{X}_t \subset X_t$.\end{proof}

\begin{theorem}
There  is  no  duality  gap  between  primal  LP  problem  \eqref{eqn:Pa}  on  measures  and
dual LP problem \eqref{eqn:Da}  on functions in the sense that $p^T = d^T$.
\label{thm:nogap}
\end{theorem}
\begin{proof}
{
According to a standard result \cite[Chapter IV, Theorem (7.2)]{barvinok} of infinite-dimensional linear programming, the zero duality gap follows from the closedness of the cone $K:=\{(\mathcal{A}\phi, \langle \phi,c \rangle_1) : \phi \in E_1^+ \}$ in $E_2 \oplus \mathbb{R}$. In order to prove the closedness of $K$, we take a sequence $\phi_k = (\mu_{0,k}, \hat{\mu}_{0,k}, \mu_k, \mu_{\partial,k})$ such that $(\mathcal{A}\phi_k, \langle \phi_k,c \rangle_1)$ weakly-* converges in $E_2 \oplus \mathbb{R}$ and show that its limit is in $K$. Since $(\mathcal{A}\phi_k)_k$ converges, $\left\langle \mathcal{A}\phi_k,\tiny{\begin{pmatrix}0 \\ 0 \\1 \end{pmatrix}} \right\rangle_2 = \mu_{0,k} + \hat{\mu}_{0,k}$ is convergent and thus bounded. This implies - using non-negativity of considered measures - that both sequences $\mu_{0,k}(X)$ and $\hat{\mu}_{0,k}(X)$ are bounded. Boundedness of $\left\langle \mathcal{A}\phi_k,\tiny{\begin{pmatrix}1 \\ 0 \\0 \end{pmatrix}} \right\rangle_2 = \mu_k(X)$ implies that $\mu_k(X)$ is bounded as well. Finally, we can remark that $\text{div}(f\mu_k)(X) = 0$ for all $k \in \mathbb{N}$, this is why boundedness of $\left\langle \mathcal{A}\phi_k,\tiny{\begin{pmatrix}0 \\ 1 \\0 \end{pmatrix}} \right\rangle_2 = $ div$(f\mu_k) + \mu_{\partial,k} - \mu_{0,k}$ and $\mu_{0,k}(X)$ implies boundedness of $\mu_{\partial,k}(X)$. Thus the four sequences of measures $\mu_{0,k}, \hat{\mu}_{0,k}, \mu_k, \mu_{\partial,k}$ are bounded in the sense of the weak-* topology. Hence, from the weak-* compactness of the unit ball (Alaoglu’s Theorem \cite[Chapter III, Theorem (2.9)]{barvinok}), there exists a subsequence $\phi_{k_i}$ that converges weakly-* to an element $\phi \in E_2$ so that $\lim\limits_{k \to \infty} (\mathcal{A}\phi_k, \langle \phi_k,c \rangle_1) = (\mathcal{A}\phi, \langle \phi,c \rangle_1)$ by weak-* continuity of $\mathcal{A}$ and $\langle .,. \rangle_1$. This proves that $K$ is closed.}
\end{proof}

\subsection{LMI approximations}

{In what follows, $\mathbb{R}_k[x]$ denotes the vector space of real multivariate polynomials of total degree less than or equal to $k$, and $\Sigma_k[x]$ denotes the cone of sums of squares (SOS) of polynomials of degree less than or equal to $k$.

Let $\kappa := \left\lceil\sum_{i=1}^{n_X} d_i / 2\right\rceil$. For $i = 0,\dots,n_X$ let $k_i := \lceil d_i / 2 \rceil$. Let $k_{min} := \max\left(k_0, \kappa\right)$ and $k \geq k_{min}$.}

Throughout the rest of this section we make the following standard standing assumption:
\begin{assumption}
One of the polynomials modeling the set $X$ is equal to $g_i(x) = R^2 - \vert x\vert^2$.
\label{assum:ball}
\end{assumption}

This assumption is completely without loss of generality since a redundant ball constraint can be always added to the description of the bounded set X.

{Problem \eqref{eqn:Da} admits a 
SOS tightening which can be written as follows:}

\begin{equation}
\begin{array}{rcl}
d_k^T \:= & \inf & \mathbf{w}'\mathbf{l} + u \: T \: \mathbf{l}_0 \\
& \text{s.t.} & u - \nabla v \cdot f = p_{0} + \sum_i p_i\: g_i \\
& & w - v - 1 = q_{0} + \sum_i q_i \: g_i \\
& & w = s_0 + \sum_i s_i \: g_i \\
& & v = t_0 + \sum_i t_i \: g_i + t \: g_1\cdots g_{n_X}
\end{array}
\label{eqn:Dak}
\end{equation}
{where the infimum is with respect to $u \geq 0$, $v, w \in \mathbb{R}_{2k}[x]$, $q_0, s_0, t_0 \in \Sigma_{2k}[x]$, $q_i,s_i,t_i \in \Sigma_{2(k-k_i)}[x]$, $i=1,\ldots,n_X$, $t \in \mathbb{R}_{2(k-\kappa)}[x]$ and $p_0, \dots, p_{n_X}$ SOS polynomials with appropriate degree. Vector $\mathbf{l}$ denotes the Lebesgue moments over $X$ indexed in the same basis in which the polynomial $w$ with vector of coefficients $\mathbf{w}$ is expressed.}

SOS problem \eqref{eqn:Dak} is a tightening of problem \eqref{eqn:Da} in the sense that any feasible solution in \eqref{eqn:Dak} gives a triplet $(u,v,w)$ feasible in \eqref{eqn:Da}.

\begin{theorem}\label{lemme8}
Problem \eqref{eqn:Dak} is an LMI problem and any feasible solution $(u_k,v_k,w_k)$ gives inner approximations $\hat{X}_t^k := \{x \in \text{int}(X) , v_k(x) + u_kt < 0\}$ of the $X_t$s. In particular, if $u_k = 0$, $\hat{X}_\infty^k := \{x \in \text{int}(X) , v_k(x) < 0\}$ is an inner approximation of $X_\infty$.
\end{theorem}
\begin{proof}
{As any linear optimization problem on SOS polynomials, \eqref{eqn:Dak} can be written as an LMI, see e.g. \cite{l10} and references therein. Any constraint of the form $\sigma = \sum\limits_{\alpha \in \mathbb{N}^n_{2r}} \sigma_\alpha x^\alpha \in \Sigma_{2r}[x]$ is indeed equivalent to the existence of a positive semi-definite matrix $Q = (q_{\beta,\gamma})_{\beta,\gamma \in \mathbb{N}^n_{r}}$ such that for any $\alpha \in \mathbb{N}^n_{2r}, \sigma_\alpha = \sum_{\beta+\gamma = \alpha}\ q_{\beta \gamma}$. The inner approximation result is a direct consequence of Lemmas \ref{lem_pos_inv} and \ref{lem_inner_approx}.}
\end{proof}

{This SOS tightening is a finite dimension convex optimization, and as such it admits a primal formulation derived from Lagrangian theory, which can be seen as an LMI relaxation of infinite dimensional LP \eqref{eqn:Pa} (see \cite[Chapter 3]{l10} for details).}

\subsection{Convergence of the inner approximations}

\begin{theorem}
Let $T > \overline{\tau}$. Then,
\begin{itemize}
    \item[1.] The sequence $(d^T_k)$ is monotonically decreasing and converging to $\lambda(X_\infty^c)$
\end{itemize}
For every $k \geq k_{min}$, let $\psi_k := (u_k,v_k,w_k)$ denote a $\frac{1}{k}$-optimal solution of the dual tightening of order $k$. One has then :
\begin{itemize}
    \item[2.] $u_k \underset{k \to \infty}{\longrightarrow} 0$
    \item[3.] $ w_k \underset{k \to \infty}{\stackrel{L^1(X)}{\longrightarrow}} \mathbb{1}_{X_\infty^c} $
\end{itemize}
\end{theorem}
\begin{proof}
\begin{itemize}
    \item[1.] {$d^T_k$ is a decreasing sequence since the sequence of feasible sets of tightening \eqref{eqn:Dak} is increasing in the sense of the inclusion (we are looking for solutions of increasing degree). \\
We are going to prove now that $\lim\limits_{k \to \infty} d^T_k =  d^T$. Let $\epsilon > 0$. Let $(u,v,w)$ be a strictly feasible triplet for problem \eqref{eqn:Da} such that $d^T \leq \int_X (w(x) + u \: T) \: d\lambda(x) \leq d^T + \frac{\epsilon}{2}$.
Such a triplet exists since problem \eqref{eqn:Da} has strictly feasible points (such as $(u^\ast,v^\ast,w^\ast) = (1,x \mapsto 1, x \mapsto 3)$).
Let $\epsilon_1 > 0$ such that:
\begin{itemize}
\item $u - \nabla v \cdot f(x) > \epsilon_1, \forall x \in X$,
\item $w(x) - v(x) - 1 > \epsilon_1, \forall x \in X$,
\item $w(x) > \epsilon_1, \forall x \in X$,
\item $v(x) > \epsilon_1, \forall x \in \partial X$.
\end{itemize}
Such an $\epsilon_1$ exists by strict feasibility of $(u,v,w)$. Let $\epsilon_2 := \min\left\{\frac{\epsilon_1}{2} ; \frac{\epsilon_1}{1 + \Vert f \Vert_{L^\infty(X)}} ; \frac{\epsilon}{2(1+T)\lambda(X)}\right\}$.
According to an extension of the Stone-Weierstrass theorem found in \cite{swh}, there exists a triplet $(\hat{u},\hat{v},\hat{w})\in \mathbb{R}_+ \times \mathbb{R}[x]^2$ such that $\vert \hat{u} - u\vert + \Vert \hat{v} - v \Vert_{L^\infty(X)} + \Vert \nabla (\hat{v} - v) \Vert_{L^\infty(X)} + \Vert \hat{w} - w \Vert_{L^\infty(X)} < \epsilon_2$.

Let $x \in X$, $\bar{x} \in \partial X$. Then, one has
$\hat{u} - \nabla\hat{v}\cdot f(x) \:  = u - \nabla v\cdot f(x) + (\hat{u} - u) - \nabla(\hat{v} - v) \cdot f(x)  > \epsilon_1 - \epsilon_2(1 + \Vert f \Vert_{L^\infty(X)}) \geq 0$,
$\hat{w}(x) - \hat{v(x)} - 1 \: = w(x) - v(x) - 1 + (\hat{w} - w)(x) - (\hat{v} - v)(x) > \epsilon_1 - 2\epsilon_2 \geq 0$,
$\hat{w}(x) = w(x) + (\hat{w} - w)(x) > \epsilon_1 - \epsilon_2 \geq 0.$
$\hat{v}(\bar{x}) = v(\bar{x}) + (\hat{v} - v)(\bar{x}) > \epsilon_1 - \epsilon_2 \geq 0$.

Hence, Assumption \ref{assum:ball} enables to use Putinar's Positivestellensatz \cite[Theorem 2.14]{l10}, giving the existence of $k \in \mathbb{N}$ and of $q_0, s_0, t_0 \in \Sigma_{2k}[x]$, $q_i,s_i \in \Sigma_{2(k-k_i)}[x]$, $i=1,\ldots,n_X$, $t_1 \in \mathbb{R}_{2(k-\kappa)}[x]$ and $p_0, \dots, p_{n_X}$ SOS polynomials with appropriate degree such that $\hat{u} - \nabla \hat{v} \cdot f = q_{0} + \sum_i q_i\: g_i$, $ \hat{w} - \hat{v} - 1 = p_{0} + \sum_i p_i \: g_i$, $ \hat{w} = s_0 + \sum_i s_i \: g_i$ and $ \hat{v} = t_0 + t_1 g_1 \cdots g_{n_X}$. This gives a feasible solution of tightening \eqref{eqn:Dak} of order $k$. Moreover, we can check that 
$\int_X(\hat{w}(y) + \hat{u} \: T) \: d\lambda(y) = \int_X(w(y)+u \: T) \: d\lambda(y) + \int_X(\hat{w} - w)(y) \: d\lambda(y) + T\:\lambda(X)\:(\hat{u}-u) \leq d^T + \frac{\epsilon}{2} + \lambda(X)\epsilon_2(1+T) \leq d^T + \epsilon$
which implies that $d^T\leq d^T_k \leq d^T+\epsilon$.
\item[2.] We define $\phi^*:=(\mu_0,\hat{\mu}_0,\mu,\mu_\partial)$ feasible for \eqref{eqn:Pa} as in the proof of Lemma \ref{LemmaInf}. In particular, $-\mathcal{A}\phi^\ast + b = (T\:\lambda(X) - \mu(X),0,0) \in K_2$ and $\langle -\mathcal{A}\phi^\ast , \psi_k \rangle + \langle b, \psi_k\rangle = \langle -\mathcal{A}\phi^\ast + b, \psi_k \rangle = (T\:\lambda(X) - \mu(X))u_k \geq 0$. On the other hand, $\langle -\mathcal{A} \phi^\ast, \psi_k \rangle = - \langle \phi^\ast, \mathcal{A}'\psi_k \rangle \leq - \langle \phi^\ast,c \rangle = -d^T$ since $\mathcal{A}'\psi_k - c \in F_1^+$ and $\phi^\ast \in E_1^+$ is optimal. Thus,
$$ 0 \leq (T\:\lambda(X) - \mu(X))u_k \leq \langle b , \psi_k \rangle - d^T. $$
According to point 1, we have that $\langle b , \psi_k \rangle \underset{k\to\infty}{\longrightarrow} d^T$, so $(T\:\lambda(X) - \mu(X))u_k \underset{k\to\infty}{\longrightarrow} 0$. Since by assumption $T > \bar{\tau} = \frac{\mu(X)}{\lambda(X)}$, this means that $u_k \underset{k\to\infty}{\longrightarrow} 0$.
}
    \item[3.] Let $\epsilon >0$. Let $t >0$ such that $\lambda(X_t \setminus X_\infty) \leq \epsilon$. Let $\Bar{k} \geq k_{min}$ such that for all $k \geq \Bar{k}$ one has that $\Vert u_k t \Vert_{L^1(X)} \leq \epsilon$ and $\vert\int_X w_k d \lambda - \lambda(X^c_\infty) \vert \leq \epsilon$. Such an integer exists from points $1$ and $2$. Using the triangular inequality and the fact that $\Vert u_k t \Vert_{L^1(X)} \leq \epsilon$ one has 
    \begin{equation}
    \Vert w_k - \mathbb{1}_{X_\infty^c}\Vert_{L^1(X)}\leq \Vert w_k + u_k t - \mathbb{1}_{X_\infty^c}\Vert_{L^1(X)} + \epsilon.
    \label{ineq}
    \end{equation}
    With the notation $\Delta = \Vert w_k + u_k t - \mathbb{1}_{X_\infty^c}\Vert_{L^1(X)}$, one has that $ \Delta = \int_{X^c_t} |w_k + u_k t - \mathbb{1}_{X_\infty^c}|d\lambda + \int_{X_t} \vert w_k + u_k t - \mathbb{1}_{X_\infty^c}\vert d\lambda. $ We denote by $\Delta_1$ and $\Delta_2$ these two terms, respectively. Using that $X_t^c \subset X_\infty^c$ and that $w_k(x) + u_k t \geq 1 + v_k(x) + u_k t \geq 1 , \forall x \in X_t^c$ (from Theorem \ref{lemme8}) we have then that
    $\Delta_1 = \int_{X^c_t} w_k + u_k t - 1 d\lambda = \int_{X^c_t} w_k d\lambda - \lambda(X^c_t) + \lambda(X^c_t)  u_k t $
    and since  $\lambda(X^c_t) u_k t \leq \Vert u_k t \Vert_{L^1(X)} \leq \epsilon$,
    \begin{equation}
        \Delta_1 \leq \int_{X^c_t} w_k d\lambda - \lambda(X^c_t) + \epsilon.
        \label{Delta1}
    \end{equation}
    Moreover, we have that
    $ \Delta_2 \leq \int_{X_t} \vert w_k \vert + \vert u_k t \vert + \vert \mathbb{1}_{X_\infty^c}\vert d\lambda$
    and therefore, using that $w_k \geq 0$ and that $\Vert u_k t \Vert_{L^1(X)} \leq \epsilon$, it holds
    $ \Delta_2 \leq \int_{X_t} w_k d\lambda + \epsilon  + \lambda(X_t \setminus X_\infty).$
    Since we have $\lambda(X_t \setminus X_\infty) \leq \epsilon$ by choice of $t$, we deduce that $ \Delta_2 \leq \int_{X_t} w_k d\lambda + 2\epsilon$. Combining this inequality with \eqref{Delta1}, we have :
    $\Delta = \Delta_1 + \Delta_2 \leq \int_X w_k d \lambda - \lambda(X^c_t) + 3 \epsilon $
    from which we deduce that $\Delta \leq 5 \epsilon$, using that $\vert\int_X w_k d \lambda - \lambda(X^c_\infty) \vert \leq \epsilon$ and  $\lambda(X_\infty^c \setminus X_t^c)\leq \epsilon$. Combining this with \eqref{ineq}, we have that  $\Vert w_k - \mathbb{1}_{X_\infty^c}\Vert_{L^1(X)}\leq 6 \epsilon $.
\end{itemize}
\end{proof}

{
\begin{remark} Despite this convergence result, one should be aware of the fact that the computational burden increases sharply with the dimension of the state space and the degree of the relaxations. Indeed, the involved polynomials have ${n+d \choose d} = {n+d \choose n}$ coefficients. Consequently, high values of $n$ and $d$ might be intractable. A possible way to handle this consists in exploiting the structure of the considered problems, such as sparsity. The key is to split the state space into low dimensional subspaces and distribute the problem over the obtained partitioning (see \cite{twlh19} as a first example of what can be done in practice for volume computation).
\end{remark}
}

\section{Numerical example}\label{sec:Numerics}

For this paper, we chose to focus on the simple example of the Van der Pol oscillator, as was done in \cite{hk14}:
\begin{equation}
\begin{cases}
\dot{x}_1 = -2 \ x_2 \\
\dot{x}_2 = 0.8 \ x_1 + 10 \ (1.02^2x_1^2 - 0.2) \ x_2.
\end{cases}
\end{equation}
Let $X = \{x \in \mathbb{R}^2 : x_1^2 + x_2^2 \leq 1\}$ and $T = \frac{100}{\pi}$.

\begin{figure}[!h]
    \centering
    \includegraphics[scale=0.7]{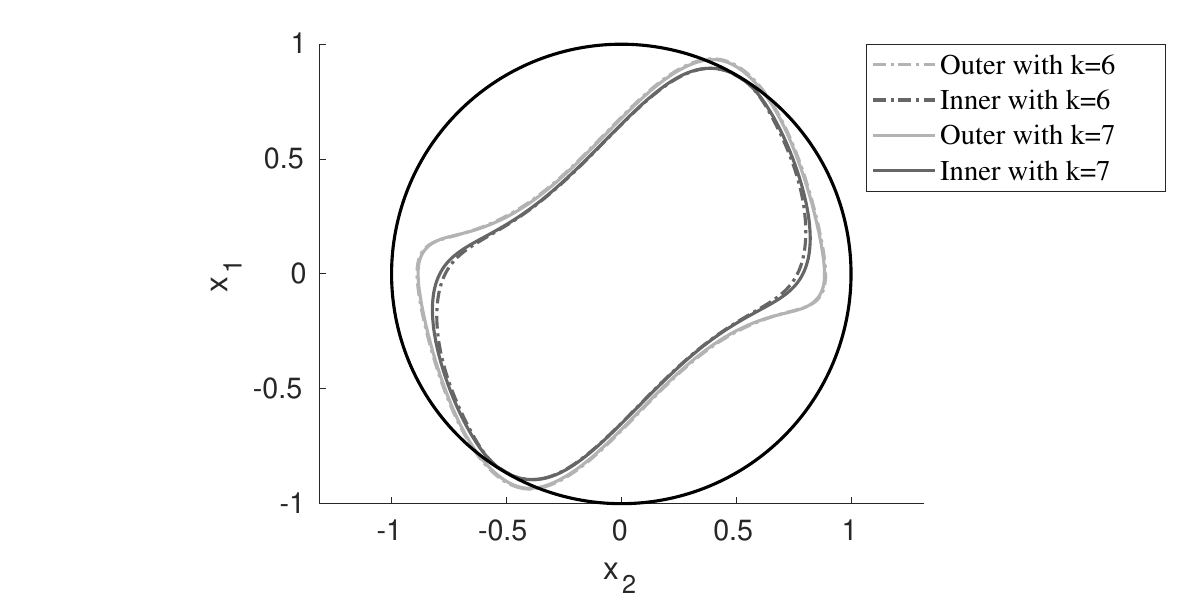}
    \caption{Outer and inner approximations of the Van der Pol MPI set in the unit disk.}
    \label{fig:MPI_VDP}
\end{figure}

We implemented the hierarchy of SOS problems \eqref{eqn:Dak} in MATLAB, using the toolbox YALMIP interfaced with the SDP solver MOSEK. For $k=6$ and $7$ (SOS degrees $12$ and $14$ respectively), we compared the obtained regions to the outer approximations computed using the framework presented in \cite{khj14}, see Figure \ref{fig:MPI_VDP}.
In this implementation, we checked at each realxation whether $u$ was near to zero: for $k = 6$, we had $u \sim 10^{-7}$, and for $k = 7$ we obtained $u \sim 10^{-6}$, which is satisfactory. {Moreover, we also ran the hierarchy with constraint $u_k = 0$ (to enforce inner approximations) and obtained the same results.}

However, we observed some difficulties:
\begin{itemize}
\item For low degrees, the only solution $v$ found by the solver is very close to the zero polynomial: the coefficients are of the order $10^{-5}$, therefore the plots are irrelevant; one loses conservativeness and several constraints are violated (namely the positivity constraint on $v$ on $\partial X$).
\item For higher degrees, the basis of monomials is not adapted since for example in dimension 1 $x^\alpha$ is close to the indicator of $\{-1,1\}$. As a result, the coefficients are of the order $10^5$ or more, and again the plots make little sense.
\end{itemize}
One can also find numerical applications of this method to actual eletrical engineering problems in \cite{ocph18} with very promising results.

\section{CONCLUSIONS}\label{sec:Conclusions}

The original motivation behind our current work is the study of transient phenomena in large-scale electrical power systems, see \cite{jmts18} and references therein. Our objective is to design a hierarchy of approximations of the MPI set for large-scale systems described by non-linear differential equations. A first step towards non-polynomial dynamics can be found in \cite{ocph18}. Since the initial work \cite{hk14} relied on the mathematical technology behind the approximation of the volume of semi-algebraic sets, we already studied in \cite{twlh19} the problem of approximating the volume of a large-scale sparse semi-algebraic set. We are now investigating extensions of the techniques for approximating the MPI set of large-scale sparse dynamical systems, and the current paper contributes to a better understanding of its inner approximations, in the small-scale non-sparse case. Our next step consists of combining the ideas of \cite{twlh19} with those of the current paper, so as to design a Lasserre hierarchy of inner approximations of the MPI set in the large-scale case, and apply it to electrical power system models.

\addtolength{\textheight}{-12cm}   




\section*{ACKNOWLEDGMENT}

We would like to thank Victor Magron and Milan Korda (LAAS-CNRS) as well as Carmen Cardozo and Patrick Panciatici (RTE) for fruitful discussions.

\end{document}